%%%%%%%%%%%%%%%%%%%%%%%%%%%%%%%%%%%%%%%%%%%%%%%%%%%%%%%%%%
%%
%%
%%     Paper  Colli--Sprekels 13
%%     On the optimal control of the viscous Cahn--Hilliard equation
%%     with hyperbolic relaxation of the chemical potential
%%
%%%%%%%%%%%%%%%%%%%%%%%%%%%%%%%%%%%%%%%%%%%%%%%%%%%%%%%%%%
%
%%%%%%%%%%%%%%%%%%%%%%%%%%%%%%%%%
%
%		COLORS FOR CORRECTIONS
%   EXAMPLE: {\pier{I want this text to become cyan}}
%   \let\pier\relax schaltet die Farbe wieder aus  

\def\pier#1{{\color{red}#1}}
\def\revis#1{{\color{red}#1}}

\def\pier#1{#1}
\def\revis#1{#1}

%\def #1{#1}
%\let\juerg\relax

%%%%%%%%%%%%%%%%%%%%%%%%%%%%%%%

\def\LaTeX{%
  \let\Begin\begin
  \let\End\end
  \def\Bcenter{\Begin{center}}
  \def\Ecenter{\End{center}}
  \let\Label\label
  \let\salta\relax
  \let\finqui\relax
  \let\futuro\relax}

\def\UK{\def\our{our}\let\sz s}
\def\USA{\def\our{or}\let\sz z}

%%%%%%%%%%%%%%%%%%%%%%%%%%%%%%%%

% scegliere fra \TeX e \LaTeX  e fra  \UK oppure \USA

%\TeX
\LaTeX

%\UK
\USA

%%%%%%%%%%%%%%%%%%%%%%%%%%%%%%%%%
%% page layout
%%%%%%%%%%%%%%%%%%%%%%%%%%%%%%%%%

\salta

\documentclass[twoside,12pt]{article}
\setlength{\textheight}{24cm}
\setlength{\textwidth}{16cm}
\setlength{\oddsidemargin}{2mm}
\setlength{\evensidemargin}{2mm}
\setlength{\topmargin}{-15mm}
\parskip2mm

%%%%%%%%%%%%%%%%%%%%%%%%%%%%%%%%%
%% packages
%%%%%%%%%%%%%%%%%%%%%%%%%%%%%%%%%

%
%    ADDED BY PIER
%
\usepackage{cite}

\usepackage{color}
\usepackage{amsmath}
\usepackage{amsthm}
\usepackage{amssymb}

\usepackage{amsfonts}
\usepackage{mathrsfs}

\usepackage{hyperref}
\usepackage[mathcal]{euscript}

\usepackage[ulem=normalem,draft]{changes}

%%%%%%%%%%%%%%%%%%%%%%%%%%%%%%%%%
%% bibliographystyle
%%%%%%%%%%%%%%%%%%%%%%%%%%%%%%%%%

%\bibliographystyle{plain}

%%%%%%%%%%%%%%%%%%%%%%%%%%%%%%%%%
%% environments
%%%%%%%%%%%%%%%%%%%%%%%%%%%%%%%%%

%

\finqui

\def\Beq{\Begin{equation}}
\def\Eeq{\End{equation}}

\def\Bthm{\Begin{theorem}}
\def\Ethm{\End{theorem}}
\def\Blem{\Begin{lemma}}
\def\Elem{\End{lemma}}

\def\Brem{\Begin{remark}\rm}
\def\Erem{\End{remark}}

\def\Bdim{\Begin{proof}}
\def\Edim{\End{proof}}
\let\non\nonumber

%%%%%%%%%%%%%%%%%%%%%%%%%%%%%%%%%
%% macros
%%%%%%%%%%%%%%%%%%%%%%%%%%%%%%%%%

% macro salvate

% sottosezioni non numerate

\def\step #1 \par{\medskip\noindent{\bf #1.}\quad}

% abbreviazioni di parole

%\def\sc{so-called}

\def\rhs{right-hand side}

% versioni inglesi (UK) o americane (USA)

% bold, cal e mathop

\def\multibold #1{\def\arg{#1}%
  \ifx\arg\pto \let\next\relax
  \else
  \def\next{\expandafter
    \def\csname #1#1#1\endcsname{{\bf #1}}%
    \multibold}%
  \fi \next}

\def\pto{.}

\def\multical #1{\def\arg{#1}%
  \ifx\arg\pto \let\next\relax
  \else
  \def\next{\expandafter
    \def\csname cal#1\endcsname{{\cal #1}}%
    \multical}%
  \fi \next}

% operatori

\def\multimathop #1 {\def\arg{#1}%
  \ifx\arg\pto \let\next\relax
  \else
  \def\next{\expandafter
    \def\csname #1\endcsname{\mathop{\rm #1}\nolimits}%
    \multimathop}%
  \fi \next}

\multibold
qwertyuiopasdfghjklzxcvbnmQWERTYUIOPASDFGHJKLZXCVBNM.

\multical
QWERTYUIOPASDFGHJKLZXCVBNM.

\multimathop
dist div dom meas sign supp .

% accorpamenti di formule citate:
% uso  \accorpa {prima}{seconda}
%      \Accorpa\cs prima seconda (con il comodo blank anche dopo)
% NB: \Accorpa definisce \cs come l'accorpamento delle due citazioni
% e scrive sul file.log

\def\Accorpa #1#2 #3 {\gdef #1{\eqref{#2}--\eqref{#3}}%
  \wlog{}\wlog{\string #1 -> #2 - #3}\wlog{}}

% macro comode

\def\<#1>{\mathopen\langle #1\mathclose\rangle}
\def\norma #1{\mathopen \| #1\mathclose \|}

%\textsl{}

\def\iot {\int_0^t}

\def\iO{\int_\Omega}

\def\iQ{\iint_Q}

\def\dt{\partial_t}
\def\dtt{\partial_{tt}}
\def\dn{\partial_{\bf n}}

\def\cpto{\,\cdot\,}

\def\checkmmode #1{\relax\ifmmode\hbox{#1}\else{#1}\fi}

% insiemi numerici

\def\erre{{\mathbb{R}}}
\def\enne{{\mathbb{N}}}

% spazi di funzioni a valori vettoriali su [0,T], [0,t], [0,s], [0,+\infty), [\delta,T]

% Come ricordare: in generale i simboli L H W  C da soli per gli spazi su (0,T)
% gli stessi raddoppiati per (0,+\infty)
% aggiunta di t o s al simbolo per (0,t) e (0,s)
% aggiunta di d al simbolo semplice o doppio per intervalli (\delta,T) e (\delta,+\infty)
% il simbolo C e i suoi derivati mettono le quadre anziche' le tonde

% Esempi   \L2V   \L\infty{V^*}   \W{1,1}H   \C0H   \LL2V   \CC0{V^*}   \Ld2V  \CCdH

\def\genspazio #1#2#3#4#5{#1^{#2}(#5,#4;#3)}
\def\spazio #1#2#3{\genspazio {#1}{#2}{#3}T0}

\def\L {\spazio L}
\def\H {\spazio H}
\def\W {\spazio W}

\def\C #1#2{C^{#1}([0,T];#2)}

% spazi di funzioni su \Omega e \Gamma

\def\Lx #1{L^{#1}(\Omega)}
\def\Hx #1{H^{#1}(\Omega)}

\def\Luno{\Lx 1}
\def\Ldue{\Lx 2}
\def\Linfty{\Lx\infty}

\def\Huno{\Hx 1}
\def\Hdue{\Hx 2}

\def\LiQ{L^\infty(Q)}

% spazi di funzioni su Q

% lettere greche

\let\theta\vartheta

\let\phi\varphi

\let\TeXchi\chi                         % new \chi, exactly on the baseline
\newbox\chibox
\setbox0 \hbox{\mathsurround0pt $\TeXchi$}
\setbox\chibox \hbox{\raise\dp0 \box 0 }
\def\chi{\copy\chibox}

% abbreviazioni specifiche del lavoro

%%%%%%%%

\def\CS{{\cal S}}
\def\CP{{\bf {(CP)}}}

\def\CU{{\cal U}}

\def\Uad{{\cal U}_{\rm ad}}

\def\CJ{{\cal J}}

\def\us{u^\alpha}
\def\phis{\phi^\alpha}
\def\mus{\mu^\alpha}
\def\ws{w^\alpha}
\def\qs{q^\alpha}
\def\ps{p^\alpha}
\def\rs{r^\alpha}

\def\VD{V^*}

\def\phiz{\phi_0}

% scelte che si possono cambiare

%\let\tau \taudef
\normalfont

%%%%%%%%%%%%%%%%%%%%%%%%%%%%%%
\Begin{document}
%%%%%%%%%%%%%%%%%%%%%%%%%%%%%%%%%

%% front page
%%%%%%%%%%%%%%%%%%%%%%%%%%%%%%%%%

\title{{\bf On the optimal control of viscous Cahn--Hilliard systems
with hyperbolic relaxation of the chemical potential}}
\author{}
\date{}
\maketitle

\Bcenter
\vskip-1.3cm
{\large\bf Pierluigi Colli$^{(1)}$}\\
{\normalsize e-mail: {\tt pierluigi.colli@unipv.it}}\\[.4cm]
{\large\bf J\"urgen Sprekels$^{(2)}$}\\
{\normalsize e-mail: {\tt juergen.sprekels@wias-berlin.de}}\\[.6cm]

%$^{(0)}$
%?????\\[.2cm]
$^{(1)}$
{\small Dipartimento di Matematica ``F. Casorati'', Universit\`a di Pavia}\\
{\small and Research Associate at the IMATI -- C.N.R. Pavia}\\ 
{\small via Ferrata 5, 27100 Pavia, Italy}\\[.2cm]
$^{(2)}$
{\small Weierstrass Institute for Applied Analysis and Stochastics}\\
{\small Mohrenstra\ss e 39, 10117 Berlin, Germany}\\[.8cm]

\Ecenter

{
\Begin{abstract}\noindent
In this paper, we study an optimal control problem for a viscous Cahn--Hilliard system with zero Neumann boundary
conditions in which a hyperbolic relaxation term involving the 
second time derivative of the chemical potential has been added to the first equation
of the system. For the initial-boundary value problem of this system, results concerning well-posedness, continuous dependence
and regularity are known. We show Fr\'echet differentiability of the associated control-to-state operator, study the associated
adjoint state system, and derive first-order necessary optimality conditions. \pier{Concerning the nonlinearities driving the system, we can include the case of logarithmic potentials. In addition, we perform an 
asymptotic analysis of the optimal control problem as the relaxation coefficient approaches zero.}
\\[2mm]
\normalfont
\\[2mm]

\End{abstract}
}
\salta

\pagestyle{myheadings}
\newcommand\testopari{\sc Colli \ --- \ Sprekels }
\newcommand\testodispari{\sc  Control \pier{for} Cahn--Hilliard systems with hyperbolic relaxation}
\markboth{\testopari}{\testodispari} 

%\markright{}

\finqui

%%%%%%%%%%%%%%%%%%%%%%%%%%%%%%%%%
%% very beginning
%%%%%%%%%%%%%%%%%%%%%%%%%%%%%%%%%

\section{Introduction}
\label{Intro}
\setcounter{equation}{0}

Let $\Omega\subset \erre^3$ be a bounded and connected domain with smooth boundary  $ \partial\Omega$. 
%and $T$ denote some final time. 
We denote by ${\bf n}$ the unit outward 
normal to $ \partial\Omega$, with the associated outward normal derivative $\,\dn$. Moreover, let $T>0$ \pier{stand for} some final time, and set
\begin{align*}
&Q_t:=\Omega \times (0,t),\quad \Sigma_t:=\partial\Omega\times (0,t), \quad \mbox{for }\,t\in(0,T],\quad\mbox{and}\quad
Q:=Q_T,\quad \Sigma:=\Sigma_T.
\end{align*}
We then study the following optimal control problem:

\vspace*{2mm}\noindent
{\bf (CP)} \,\,Minimize the tracking-type cost functional
\begin{align}
\label{cost}
{\cal J}(\phi, u)\,&:=\,\frac{b_1}2\iQ|\phi -\phi_Q|^2\,+\,\frac{b_2}2 \iO|\phi(T)-\phi_{\Omega}|^2\,+\,\frac{b_3}2\iQ|u|^2
\,+\,\kappa\,G(u)\nonumber\\
&=:J(\phi,u)\,+\,\kappa\,G(u)
\end{align}
subject to the initial-boundary value system 
\begin{align}
\label{ss1}
&\alpha \dtt\mu + \dt\phi - \Delta\mu=0 &&\mbox{a.e. in }\,Q, \\
\label{ss2}
&\tau\dt\phi-\Delta \phi + f'(\phi) = \mu + w  &&\mbox{a.e. in }\,Q,\\
\label{ss3}
&\gamma \dt w + w = u   &&\mbox{a.e. in }\,Q,\\
\label{ss4}
&\dn\mu=\dn\phi = 0 &&\mbox{a.e. on }\,\Sigma,\\
\label{ss5}
&\mu(0)=\mu_0,\ \,  (\dt \mu) (0)=\nu_0,\ \,\phi(0)=\phi_0,\ \,w(0)= w_0 &&\mbox{a.e. in }\,\Omega,
\end{align}
\Accorpa\State ss1 ss5
and to the control constraint
\begin{equation}
\label{defUad}
u\in \Uad=\{u\in\CU: \ \underline u(x,t)\le u(x,t)\le \overline u(x,t) \,\mbox{ for a.a. $(x,t)$\, in }\,Q\}.
\end{equation}
Here, the control space is specified by
\begin{equation}
\label{defU}
\CU=\LiQ.
\end{equation}
The given bounds $\,\underline u,\overline u\in {\LiQ}\,$ satisfy $\,\underline u\le\overline u$ 
almost everywhere in $Q$, the targets $\phi_Q, \, \phi_{\Omega}$ are given functions, and $b_1\ge 0$, $b_2\ge 0$, $b_3>0$ are constants.
Moreover, $G:L^2(Q)\to [0,+\infty)\,$ is a convex and continuous functional enhancing the occurrence of \emph{sparsity} of optimal controls, 
where $\kappa\ge 0\,$ is a fixed constant (the sparsity parameter); a typical choice for $G$ is
\begin{equation}
\label{defG}
G(u)=\|u\|_{L^1(Q)}=\iint_Q |u|\,.
\end{equation}

The equations \eqref{ss1}--\eqref{ss2} constitute a variation of the well-studied  viscous Cahn--Hilliard system 
(see, e.g., the recent contributions 
\cite{CoSpTr, CoGaMi, GiRoSi,GSS1,KO} and the references given therein)
\begin{align}
\label{ch1}
&\dt\phi - \Delta \mu = 0 \quad \mbox{in }\,Q, \\
\label{ch2}
& \tau\dt\phi-\Delta \phi + f'(\phi) = \mu + w  \quad \mbox{in }\,Q,
\end{align}
in which a hyperbolic relaxation term $\alpha \dtt\mu$  has been added in the first equation. While 
hyperbolic relaxations of the viscous Cahn--Hilliard system involving an inertial term of the phase variable $\phi$ have been considered
in the previous works \cite{Bon1,Bon,Chen,CMY,DMP,SZ,SK},  an inertial term like 
$\alpha \dtt\mu$ in \eqref{ss1} was to our knowledge studied only in the recent paper \cite{CS12}.
 
In the above class of problems (which model, e.g., the coarsening processes in binary metallic alloys), the unknown 
functions $\,\phi\,$ and $\,\mu\,$ usually stand for the \emph{order parameter}, which can represent a scaled
density of one of the involved phases, and the \emph{chemical potential}
associated with the phase separation process, respectively. 
The state variables 
$\,\phi\,$ and $\,\mu\,$  are monitored through 
the input variable $\,w$, which is in turn determined by the action of the control $\,u\,$ via the linear control equation~\eqref{ss3}. 
Eq.~\eqref{ss3} models how the ``forcing'' $w$ is generated by the external control $u$. We remark that Eq. \eqref{ss3}
could be replaced by much more general partial differential equations modeling the relation between an $L^\infty$-control $u$ and a forcing $w$.
In the practical application to the coarsening process of a binary metallic alloy, a typical aim of the control problem \CP\ is to monitor the 
system in such a way as to achieve desired distributions $\phi_Q$ and $\phi_\Omega$ of one of the metallic components in the container
$\,\Omega$ during the process and at the final time $t=T$, respectively, at minimal cost.  

In the state system \eqref{ss1}--\eqref{ss5},   $\gamma>0 $ and $\tau>0$ are given fixed constants,  and
$\mu_0,\,\,\nu_0,\,\,\phi_0,\,\,w_0$ are given initial data. Moreover,
$f' \,$ denotes the derivative 
of a double-well potential $f$ and stands for the local part of the thermodynamic force driving the evolution of the system.
Typically, $f$  is split  into a (possibly nondifferentiable) convex part~$f_1$ 
 and a smooth and concave perturbation~$f_2$. Typical and physically relevant examples for $\,f\,$ 
are the so-called {\em classical regular potential}, the {\em logarithmic double-well potential\/},
and the {\em double obstacle potential\/}, which are given, in this order,~by
\begin{align}
  & f_{\rm reg}(r) := \frac 14 \, (r^2-1)^2 \,,
  \quad r \in \erre, 
  \label{regpot}
  \\
  & f_{\rm log}(r):=\left\{\begin{array}{ll}
(1+r)\,\ln(1+r)+(1-r)\,\ln(1-r)-c_1r^2 &\quad\mbox{if }\,r\in(-1,1)\\
2\ln(2)-c_1 &\quad\mbox{if }\,r\in\{-1,1\}\\
+\infty&\quad\mbox{if }\,r\not\in [-1,1]
\end{array}, 
\right. 
  \label{logpot}
  \\[1mm]
  & f_{\rm 2obs}(r):=\left\{\begin{array}{ll}
c_2(1-r^2)  &\quad\mbox{if }\,r\in[-1,1]\\
+\infty&\quad\mbox{if }\,r\not\in [-1,1]
\end{array}. \right.
  \label{obspot}
\end{align}
Here, the constants $c_i$ in \eqref{logpot} and \eqref{obspot} satisfy
$c_1>1$ and $c_2>0$, so that $f_{\rm log}$ and $f_{\rm 2obs}$ are nonconvex.
Notice that for $f=f_{\rm log}$ the term $f'(\phi)$
occurring in \eqref{ss2} becomes singular as $\phi\searrow -1$ and $\phi\nearrow1$, which forces the order parameter $\phi$ to attain its values in the physically meaningful range~$(-1,1)$.  
In the nonsmooth case \eqref{obspot}, the convex part~$f_1$ is given by
the~indicator function of $[-1,1]$.
Accordingly, in such cases one has to replace the derivative of the convex part
by the subdifferential $\partial f_1$ and, consequently, to interpret \eqref{ss2} as a differential inclusion
or a variational inequality.
We admit in this paper only the regular and logarithmic cases.
 
There exists a large literature concerning the optimal control of viscous Cahn--Hilliard systems that cannot be cited here; 
in this connection, we refer the interested reader to the references given in the recent papers \cite{CGM, CGSCC,CoSpTr}. On the other hand, optimal control problems 
like \CP\ for the system \State\ have not been investigated before. \revis{In particular, a standout reference is paper  
\cite{CoSpTr} which deserves special recognition for its in-depth exploration of the optimal control problem applied to 
system~\State\ with $\alpha=0$, including the case of the logarithmic potential \eqref{logpot}. Notably, the resulting viscous 
Cahn--Hilliard system is far more manageable, as the absence of the inertial term $\alpha \dtt\mu$ (when $\alpha > 0 $) 
transforms the challenge from a hyperbolic system into a more tractable parabolic one. In contrast, here we take on the more 
complex scenario, successfully demonstrating convergence results, also for the optimal control problem, as the coefficient $\alpha$ 
tends to zero.}

The paper is organized as follows. In the following section, we formulate the general assumptions and state the main results concerning 
the system~\State.  In Section~3, we show the Fr\'echet differentiability of the control-to-state operator associated with the state 
system. Section 4  then brings an analysis of the control problem \CP\ for fixed $\alpha\in (0,1]$. 
Besides existence of an optimal control, first-order necessary optimality conditions are shown in terms of a variational 
inequality and the corresponding adjoint state variables. We also obtain a sparsity result in the case $\kappa>0$.
The final section  then deals with asymptotic results as $\alpha\searrow0$. We show that both the state and adjoint state variables 
converge in a well-defined sense to their counterparts associated the nonrelaxed case when $\alpha=0$, and we also prove
a convergence result for the optimal controls as $\alpha\searrow0$. 

Prior to this, let us  fix some notation.
For any Banach space $X$, we let $X^*$  denote 
its dual space, and $\norma{ \cpto }_X$ stands for the norm 
in $X$ and any power of $X$.   
For two Banach spaces $X$ and $Y$ that are both continuously embedded in some topological vector space~$Z$, the linear space
$X\cap Y$ is the Banach space equipped with its natural norm $\norma v_{X\cap Y}:=\norma v_X+\norma v_Y\,$  for $v\in X\cap Y$.
The standard Lebesgue and Sobolev spaces $L^p(\Omega)$ and $W^{m,p}(\Omega)$ are defined on $\Omega$ for
$1\le p\le\infty$ and $\,m \in \enne \cup \{0\}$. 
For the sake of 
convenience, we denote the norm of $\,L^p(\Omega)\,$ by $\,\|\,\cdot\,\|_p\,$ for $1\le p\le\infty$. 
If $p=2$, we employ the usual notation $H^m(\Omega):= W^{m,2}(\Omega)$. We also set
\begin{equation}
  H := \Ldue , \quad V := \Huno, \quad W:=\bigl\{v\in \Hdue : \ \dn v =0 \, \hbox{ on }\, \Gamma \bigr\}.
  \end{equation}
Moreover, $\VD$~is the dual space of $V$,
and $\<\cpto,\cpto>$ stands for the duality pairing between $\VD$ and~$V$. We also denote by $(\,\cdot\,,\,\cdot\,)$ the natural inner product in $\,H$.
As usual, $H$ is identified with a subspace of the dual space $\VD$ according to the identity
\begin{align*}
	\langle u,v\rangle =(u,v)
	\quad\mbox{for every $u\in H$ and $v\in V$}.
	\end{align*}
Note that $W \subset V \subset H\equiv H^* \subset \VD$ with dense and compact embeddings. 
Finally, we introduce for functions $v\in L^1(Q)$ the temporal antiderivative $\,1\ast v\,$ by setting
\begin{equation}
\label{anti}
(1\ast v)(x,t):=\int_0^t v(x,s)\,ds \quad\mbox{for a.e. $\,x\in\Omega\,$ and all $\,t\in [0,T]$}\,.
\end{equation}

About the constants used in the sequel for estimates, we adopt the rule that $\,C\,$ denotes any 
positive constant that depends only on the given data. The value of such generic constants $\,C\,$ 	may 
change from formula to formula or even within the lines of the same formula. Finally, the notation $C_\delta$ indicates a positive constant that additionally depends on the quantity $\delta$.   

%%%%%%%%%%%%%%%%%%%%%%%%%%%%%%%%%%%%%%%%%%%%%%%%%%%%%%%%%%%%%%%%%%%%%%%%

\section{General assumptions and properties of the state system}
\setcounter{equation}{0}
In this section, we formulate the general assumptions for the data of the system \State\ 
and state known existence, continuous dependence, and regularity results.
First, let us remark that the positive parameter $\alpha$ is not listed in the assumptions below, 
since it is also involved in the related asymptotic analysis. Consequently, we let 
$$ 0 < \alpha \leq 1 . $$ 
On the other hand, throughout the paper we generally assume:

\vspace*{2mm}\noindent
{\bf (A1)} \quad $\tau>0$ and $\gamma>0$ are fixed constants.

\vspace*{2mm}\noindent
{\bf (A2)} \quad $f=f_1+f_2$, where \,$f_1 :\erre\to [0,+\infty]\,$  is convex and lower semicontinuous with $f_1(0)=0$, and where
$\,f_2:\erre\to\erre$\, has a Lipschitz continuous first derivative $f_2'$ on $\erre$.

\vspace*{1mm} \noindent
Note that {\bf (A2)} implies that the subdifferential \,$\partial f_1\,$ is maximal monotone 
in $\erre \times \erre$ and satisfies \,$0 \in \partial f_1(0)$. Let $\,D(\partial f_1)\,$ denote the domain of the subdifferential.
An important requirement is the following:

\vspace*{2mm}\noindent  
{\bf (A3)} \quad There are $\,-\infty\le r_-<0<r_+\le +\infty\,$ such that $\,D(\partial f_1)=(r_-,r_+)$\, and the restrictions of
$f_1$ and $f_2$  to $(r_-,r_+)$ belong to $C^3(r_-,r_+)$. Hence, for $r\in (r_-,r_+)$, we have $\partial f_1(r)=\{f_1'(r)\}$. Moreover, we assume that 
\begin{align}
\label{singular}
\lim_{r\searrow r_-} f'_1(r)=-\infty, \quad \lim_{r\nearrow r_+} f_1'(r)=+\infty.
\end{align}

\vspace*{1mm}\noindent
Observe that the conditions {\bf (A2)} and {\bf (A3)} are fulfilled in each of the cases considered in \eqref{regpot} and \eqref{logpot}, 
with the domain $ D(\partial f_1)$ given by $\erre$ and $(-1,1)$, respectively.

We further assume: 

\vspace*{2mm} \noindent
{\bf (A4)} \quad $\mu_0\in W$, $\nu_0\in V$, \pier{$w_0\in L^\infty(\Omega)$}, $\phi_0\in W$, and it holds \,$r_-<\phi_0(x)<r_+$ \,for all $x\in\overline\Omega$.
 
\vspace*{2mm} \noindent
Observe that $\phi_0 \in  W $ implies that $\phi_0\in C^0(\overline \Omega)$. Moreover, we obviously have that
$$\,m_0 := \frac 1 {|\Omega|}\iO \phiz\,$$ lies in the interior of  \,$D(\partial f_1)$. 
Here, $\,|\Omega|\,$ denotes the Lebesgue measure of $\,\Omega$,  and $m_0$ thus represents the mean value of $\phiz$. 
In the following, we use the general notation $\overline  v$ to denote the mean value of a generic function $v\in\Luno$. If $v$ is in $V^*$,
then we can set
\Beq
  \overline v : = \frac 1 {|\Omega|} \, \langle v , 1 \rangle
  % \quad \hbox{for every $v\in V^*$},
  \label{defmean}
\Eeq
as well, noting that the constant function 1 is an element of $V$. Clearly, $\overline  v$ is the usual mean value of $v$ if $v\in H$, and 
$m_0 = \overline \phiz$. 

\pier{Finally, once and for all we} fix a ball in the control space that contains $\Uad$:

\vspace*{2mm}\noindent
{\bf (A5)} \quad The constant $R>0$ is such that $\,\Uad\subset {\cal U}_R:=\{u\in\LiQ: \|u\|_{\LiQ}<R\}$.

\vspace*{2mm} 
In the following theorem, we collect results that have been stated and proved in Theorems~2.2 to 2.5 in \cite{CS12}.
Observe that some of these results can be shown under slightly weaker assumptions, but the version stated here is tailored
to the application to the control problem \CP. We have the following result.
\Bthm
\label{due.uno}
Suppose that {\bf (A1)}--{\bf (A5)} are satisfied. Then the state system  \State\ has for any $u\in {\cal U}_R$ and every $\alpha\in (0,1]$ a unique
solution triple $(\mu^\alpha,\phi^\alpha,w^\alpha)$ satisfying the regularity requirements
\begin{align}
  & \mu^\alpha \in W^{2,\infty}(0,T;H)\cap W^{1,\infty}(0,T;V) \cap \L\infty W,
  \label{regmu}
  \\
  & \phi^\alpha \in \W{1,\infty}H \cap \H1V \cap \L\infty W,
  \label{regphi}
  \\
  & w^\alpha \in W^{1,\infty}(0,T;\pier{\Linfty}).
  \label{regw}
\end{align}
\Accorpa\Regsoluz regmu regw
Moreover, there is a constant $K_1>0$, which depends only on the data of the state system and $R$ and not on $\alpha\in (0,1]$, such that
\begin{align}
&\alpha\norma{\mu^\alpha}_{\W{2,\infty}{V^*}} + 
\alpha^{1/2}\norma{\mu^\alpha}_{\W{1,\infty}H} 
+ \norma{\mu^\alpha}_{\L\infty V} + \|w^\alpha\|_{W^{1,\infty}(0,T;\pier{\Linfty})}\non\\
&+ \norma{\phi^\alpha}_{\W{1,\infty}H \cap \H1V \cap \L\infty W}
+ \|f_1'(\phi^\alpha)\|_{\L\infty H}
  \leq K_1\,.
  \label{stimaK1} 
\end{align}
In addition, there holds a uniform separation condition in the following form: for any $\alpha\in (0,1]$, there are 
constants $r_*(\alpha),r^*(\alpha)$, which depend only on the data of the state system and $R$, such that
\begin{equation}
\label{separa}
r_-< r_*(\alpha) \le \phi^\alpha(x,t) \le r^*(\alpha) < r_+ \quad\mbox{for all }\,(x,t)\in\overline Q.
\end{equation}
\Ethm
The above result implies that for every $\alpha\in (0,1]$ the control-to-state operator
\begin{equation}
\label{defSal}
\CS^\alpha=(\CS_1^\alpha,\CS_2^\alpha,\CS_3^\alpha):u\mapsto \CS^\alpha(u)=(\CS_1^\alpha(u),\CS_2^\alpha(u),\CS_3^\alpha(u))
:=(\mu^\alpha,\phi^\alpha,w^\alpha)
\end{equation}
is a well-defined mapping between $\CU_R$ and the Banach space \pier{identified} by the regularity 
properties~\eqref{regmu}--\eqref{regw}.
\Brem
\label{PIER1}
Notice that, due to the compactness of the embedding $W\subset C^0 (\overline{\Omega}) $, it follows from \cite[Sect.~8, Cor.~4]{Simon} \pier{and \eqref{stimaK1}} that 
\pier{$\phi^\alpha \in C^0(\overline Q)$}. By the same token, the compact embedding $V\subset L^p(\Omega)$ for $1\le p<6$ yields that also 
$\mu^\alpha\in C^1([0,T];L^p(\Omega))$ for $1\le p<6$. \pier{In addition,} owing to \eqref{separa} and {\bf (A3)} there exists for every $\alpha\in(0,1]$ a constant $K_2(\alpha)>0$, which depends only
on $R$ and the data of the state system, such that 
\begin{equation}
\label{ssboundal}
\max_{i=1,2}\,\max_{j=0,1,2,3}\,\|f_i^{(j)}(\phis)\|_{\LiQ}\,\le\,K_2(\alpha)\,,
\end{equation}
whenever $\,\phis=\CS_2^\alpha(u)\,$ for some $u\in\CU_R$. \pier{Moreover, in the case when 
$\,r_- = -\infty$ and $r_+ =  +\infty\,$ so that $\,D(\partial f_1)= \erre$\, and \,$f_1 \in C^3 (\erre)$\, in {\bf (A3)}, the values $r_*(\alpha)$ and $ r^*(\alpha)$ 
in \eqref{separa}, as well as the constant $K_2(\alpha)$ in \eqref{ssboundal}, 
can be taken independent of $\alpha$. This is due to the bound \eqref{stimaK1} which ensures that
$ \norma{\phi^\alpha}_{C^0(\overline Q)} \leq C.$}
\Erem

We also have the following continuous dependence result.
\Bthm
\label{due.tre}
Suppose that {\bf (A1)}--{\bf (A5)} are fulfilled, and let $\alpha\in (0,1]$. 
Then there is a constant $K_3(\alpha)>0$ which depends only on $R$ and the data of the state system such that the following holds true:
whenever $u_i\in\CU_R$ are given and $(\mu_i^\alpha,\phi_i^\alpha,w_i^\alpha)=\CS^\alpha(u_i)$, $i=1,2$, then   
\begin{align}
&\norma{\mu_1^\alpha-\mu_2^\alpha}_{\H2\VD\cap\W{1,\infty} H \cap \L\infty V}
+ \norma{\phi_1^\alpha-\phi_2^\alpha}_{\H1H\cap\L\infty V\cap\L2W}
\non \\
&+ \|w_1^\alpha-w_2^\alpha\|_{H^1(0,T;H)}\,\leq\, K_3(\alpha) \norma{u_1-u_2}_{\L2H}\,.
  \label{stimaK4}
\end{align} 
\Ethm 
\Bdim
First observe that we have $(w_1^\alpha-w_2^\alpha)(0)=0$ and 
$$\gamma\dt(w_1^\alpha-w_2^\alpha)+w_1^\alpha-w_2^\alpha=u_1-u_2 \quad\mbox{a.e. in }\,Q.$$
Testing this identity by $\,\dt(w_1^\alpha-w_2^\alpha)\,$ then immediately yields that
$$\|w_1^\alpha-w_2^\alpha\|_{H^1(0,T;H)}\,\le\,C\,\|u_1-u_2\|_{\L2H},$$
where the constant $C$  depends on $R$ and the data of the state system, but not on $\alpha$.
The validity of \eqref{stimaK4} is then a direct consequence of the inequality (2.27) in the
statement of \cite[Thm.~2.5]{CS12}. 
\Edim
\section{Differentiability of the solution operator}
\setcounter{equation}{0}

In this section, we are going to prove differentiablity properties for the solution operator $\CS^\alpha$, where, throughout this section,
we assume that $\alpha\in (0,1]$ is fixed. To this end, let $u\in \CU_R$ be a fixed control with associated state $(\mus,\phis,\ws)
=\CS^\alpha(u)$. We then consider the linearization of the state system \State\
at $u$ in the direction $h\in \L2H$, which we write in the form
\begin{align}
\label{ls1}
&\alpha \langle\dtt\eta,\rho\rangle +\iO\dt\psi\rho+ \iO\nabla\eta\cdot\nabla\rho=0 &&\mbox{for all $\rho\in V$ and a.e. $t\in(0,T)$,}\\
\label{ls2}
&\tau\dt\psi-\Delta\psi+f''(\phis)\psi=\eta+v&&\mbox{a.e. in $Q$,}\\
\label{ls3}
&\gamma\dt v+v=h &&\mbox{a.e. in $Q$,}\\
\label{ls4}
&\dn\psi=0 &&\mbox{a.e. on $\Sigma$,}\\
\label{ls5}
&\eta(0)=\dt\eta(0)=\psi(0)=v(0)=0&&\mbox{a.e. in $\Omega$}.
\end{align}
\Accorpa\LinState ls1 ls5

\noindent Notice that \eqref{ls1} is the weak form of the PDE $\,\alpha\dtt\eta+\dt\psi-\Delta\eta=0\,$ together with 
the boundary condition $\dn\eta=0$.

We have the following result concerning well-posedness.
\Bthm
\label{tre.uno}
Suppose that {\bf (A1)}--{\bf (A5)} are fulfilled, let $\alpha\in(0,1]$ and $u\in\CU_R$ be given, as well as $(\mus,\phis,\ws)=\CS^\alpha(u)$. 
Then the system \LinState\ has for every $h\in\L2H$ a unique solution triple $(\eta_h,\psi_h,v_h)$ such that
\begin{align} 
\label{regeta}
&\eta_h\in W^{2,\infty}(0,T;V^*)\cap W^{1,\infty}(0,T;H)\cap L^\infty(0,T;V),\\
\label{regpsi}    
&\psi_h\in W^{1,\infty}(0,T;H)\cap H^1(0,T;V)\cap L^\infty(0,T;W),\\
\label{regv}
&v_h\in H^1(0,T;H). 
\end{align}
Moreover, there is a constant $K_4(\alpha)$, which depends on $R$ and the data, such that
\begin{align}
\label{contin}
&\|\eta_h\|_{W^{2,\infty}(0,T;V^*)\cap W^{1,\infty}(0,T;H)\cap L^\infty(0,T;V)}
+ \|\psi_h\|_{W^{1,\infty}(0,T;H)\cap H^1(0,T;V)\cap L^\infty(0,T;W)} \non\\
&+\|v_h\|_{H^1(0,T;H)}\,\le\,K_4(\alpha)\|h\|_{\L2H}\,. 
\end{align}
\Ethm
\Bdim 
Obviously, the PDE \eqref{ls3}, together with the initial condition $v(0)=0$ a.e. in $\Omega$, has for every $h\in\L2H$ the unique solution
$$
v_h(x,t)=\frac 1\gamma\int_0^t e^{-(t-s)/\gamma}h(x,s)\,ds\quad\mbox{for a.e. }\,(x,t)\in Q.$$
It is readily seen that
\begin{equation}
\label{bonucci}
\|v_h\|_{\H1H}\,\le\,C\|h\|_{\L2H},
\end{equation}
where, here and in the remainder of this proof, $C>0$ denotes constants that may depend on $\alpha$, $R$, and the data of the system, but neither
on $n\in\enne$ nor on $h\in\L2H$. We thus have to show the existence of a \pier{unique} pair $(\eta_h,\psi_h)$ that has zero initial data and solves  
\eqref{ls1}, \eqref{ls2} with $v=v_h$, and \eqref{ls4}. 
 
We argue by a Faedo--Galerkin approximation
using a special basis. 
To this end, we take the eigenvalues $\,\{\lambda_j\}_{j\in\enne}\,$ of the eigenvalue problem 
$$-\Delta y=\lambda y \quad\mbox{in }\,\Omega, \qquad \dn y=0 \quad\mbox{on }\,\partial \Omega,$$
and let $\,\{e_j\}_{j\in\enne}\subset W\,$ be associated eigenfunctions, normalized by $\,\|e_j\|_H=1$, $j\in\enne$. Then\pier{, it turns out that} 
\begin{align*}
&0=\lambda_1<\lambda_2\le \ldots, \qquad \lim_{j\to\infty}\lambda_j=+\infty,\\
&\iO e_je_k=\iO \nabla e_j\cdot\nabla e_k=0\quad \mbox{for $\,j\not= k$},
\end{align*}
and \pier{it is clear} that $e_1$ is just the constant function $|\Omega|^{-1/2}.$
We then define the $n$-dimensional spaces $\,V_n:={\rm span}\{e_1,\ldots,e_n\}$ for $\,n\in\enne$, where $V_1$ is just 
the space of constant functions on $\,\Omega$. 
It is well known that the union of these spaces is dense in both $\,H\,$ 
and~$\,V$.  
The approximating $n$-dimensional problem is stated as follows:
find functions
\begin{equation}
\label{discrete}
\eta_n(x,t)=\sum_{j=1}^n \eta_{nj}(t)e_j(x),\quad \psi_n(x,t)=\sum_{j=1}^n \psi_{nj}(t)e_j(x),
\end{equation}
such that 
\begin{align}
  & \alpha ( \dtt\eta_n(t) , \rho ) 
  + ( \dt\psi_n(t), \rho) 
  + \iO \nabla\eta_n(t) \cdot \nabla \rho
  = 0
  \non
  \\
  & \quad \hbox{for all $t\in [0,T]$ and every $\rho\in V_n$},
  \label{ls1n}
  \\[1mm]
  &  \tau (\dt \psi_n (t), \rho)
  + \iO \nabla\psi_n (t) \cdot \nabla \rho
  + (\pier{(f''(\phis)\psi_n)}(t), \rho) = 
  (\eta_n(t) + v_h(t), \rho) 
  \non
  \\
  & \quad \hbox{for all $t\in [0,T]$ and every $\rho\in V_n$},
  \label{ls2n}
  \\[2mm]
  &\eta_n(0)=\dt \eta_n (0)= \psi_n(0)= 0 \quad \hbox{ a.e.~in }{\Omega} \,.
  \label{ls3n}
\end{align}
We now take $\rho=e_k$ in the equations \eqref{ls1n} and \eqref{ls2n}, for $k=1,\ldots,n$, obtaining 
%the system
%\begin{align}
%\label{ODE1}
%&\alpha \frac{d^2}{dt^2} \mu_{nk} + \frac{d}{dt} 
%\phi_{nk}+\lambda_k\,\mu_{nk}=0 \quad \mbox{in }\,(0,T),\\
%\label{ODE2}
%&\tau\frac{d}{dt} \phi_{nk}+\lambda_k\,\phi_{nk}
%+(\betaeps(\phin) + f_2 (\phin),e_k)= \mu_{nk} +
%(g, e_k) 
%\quad \mbox{in }\,(0,T),\\
%\label{ODE4}c\mdseries
%& \mu_{nk}(0)=(\mu_0,e_k), \quad \frac{d}{dt}\mu_{nk}(0)=(\nu_0,e_k), \quad \phi_{nk}(0)=(\phiz,e_k).&&{}
%\end{align}
a Cauchy problem for an explicit system of linear ordinary differential equations with zero initial conditions, which is of second order
in the variables $\eta_{nk}$ 
and of first order in the variables $\psi_{nk}$. Moreover, the coefficient functions and source terms all belong to
$H^1(0,T)$. By Carath\'eodory's theorem, the Cauchy problem has a unique solution expressed by $\eta_{nk}, \, \psi_{nk}$, 
with $ \eta_{nk} \in H^3(0,T)$ and $ \psi_{nk} \in H^2(0,T)$, for $k=1,\ldots,n$. This solution 
uniquely determines a pair $(\eta_n,\psi_n)\in H^3(0,T;V_n)\times H^2(0,T;V_n)$ that solves \eqref{ls1n}--\eqref{ls3n}. 

We now derive a series of a priori
estimates for the finite-dimensional approximations. 

\step
First estimate

We take the time derivative of \eqref{ls2n} and then test it by $\dt \psi_n$. Then, we choose $\rho=\dt\eta_n$ in \eqref{ls1n}
and add the resulting equations, noting the cancellation of two terms. Next, we integrate with respect to time over $[0,t]$ for
$t\in(0,T]$, obtaining the identity 
\begin{align}
&\frac \alpha 2\|\dt \eta_n(t)\|_H^2 + \frac 1 2\|\nabla \eta_n(t)\|_H^2 +
\frac \tau 2\|\dt \psi_n(t)\|_H^2 +\iint_{Q_t} |\nabla(\dt\psi_n)|^2 \nonumber \\
&= -\iint_{Q_t} f'' (\phis)|\dt\psi_n|^2 -\iint_{Q_t}f'''(\phis)\dt\phis \,\psi_n\,\dt\psi_n+ \iint_{Q_t} \dt v_h \, \dt\psi_n\non\\
&=: I_1+I_2+I_3,
\label{help1}
\end{align}
with obvious notation. Owing to \eqref{separa}, Young's inequality, and \eqref{bonucci}, we have 
$$
I_1+I_3\,\le\,C\iint_{Q_t}|\dt\psi_n|^2 \,+\,C\|h\|_{\pier{\L2H}}^2\,.
$$
Now \pier{observe} that, thanks to the zero initial conditions \pier{and H\"older's inequality},
\begin{equation}
\label{help2}
\|\eta_n(t)\|_H^2\,\le\,t\iint_{Q_t}|\dt\eta_n|^2, \quad \|\psi_n(t)\|_H^2\,\le\,t\iint_{Q_t}|\dt\psi_n|^2, \quad\mbox{for all }\,t\in[0,T].
\end{equation}
Therefore, by virtue of \eqref{separa}, \pier{\eqref{ssboundal}}, Young's inequality, and the continuity of the embedding $V\subset L^4(\Omega)$, 
\begin{align} 
I_2&\le C\iot \|\dt\phis(s)\|_{L^4(\Omega)}\,\|\psi_n(s)\|_H\,\|\dt\psi_n(s)\|_{L^4(\Omega)}\,ds \non\\
&\le \,\frac 12\iint_{Q_t}\bigl(|\nabla\dt\psi_n|^2\,+\,|\dt\psi_n|^2\bigr)\,+\,C\iot\|\dt\phis(s)\|_V^2\,\pier{\|\dt\psi_n\|_{\L2H}^2}\,ds,
\end{align}
where, in view of \eqref{stimaK1}, the function $\,\,s\mapsto \|\dt\phis(s)\|_V^2\,\,$ belongs to $L^1(0,T)$. Hence, invoking %Young's inequality and 
\pier{Gronwall's lemma}, it is straightforward to infer 
from the above estimates that
\begin{align}
&\alpha^{1/2} \|\eta_n\|_{W^{1,\infty}(0,T; H)} + \|\nabla \eta_n\|_{{\L\infty H}^3} \non \\
&{}+
\|\pier{\psi_n}\|_{W^{1,\infty}(0,T;H)\cap \H1V} \,\leq\,C \norma{ h }_{\L2 H}\,.
\label{esti2}
\end{align}

\step 
Complementary estimates

In the following, let $P_n$ denote the $H$-orthogonal projection operator onto $V_n$. It is well known that
\begin{equation}
\label{Otto}
\|P_n\rho\|_H \,\le\,\|\rho\|_H \quad\mbox{and}\quad \|P_n\rho\|_V\,\le\,C_\Omega\|\rho\|_V\quad\mbox{for all }\,\rho\in V,
\end{equation}
with a constant $C_\Omega>0$ that depends only on $\Omega$.

Taking now an \pier{arbitrary} $\rho\in V$ in \eqref{ls1n}, and using \eqref{Otto} and the orthogonality of the basis functions, we find that 
\begin{align}
&\alpha \< \dtt\eta_n (t) ,\rho > =  
\alpha ( \dtt\eta_n (t) , P_n (\rho) )  + \alpha ( \dtt\eta_n (t) , \rho - P_n (\rho) )
\non\\
&\leq |\alpha ( \dtt\eta_n (t) , P_n (\rho) ) |
\leq \biggl|  ( \dt\psi_n (t), P_n(\rho) )  + \iO \nabla\eta_n (t) \cdot \nabla P_n(\rho)
\biggl| \non \\
&\leq C \bigl(\|\dt \psi_n\|_{\L{\infty}H} +  \|\nabla \eta_n \|_{{\L{\infty}H}^3} \bigr)  \norma{\rho}_V \quad \hbox{for a.e. } t\in (0,T), 
\label{estP1}
\end{align}
so that from \eqref{esti2} it follows that
\begin{equation}
\label{estP2}
\alpha \|\dtt \eta_n\|_{\L{\infty}{\VD}} \,\le\,C\|h\|_{\L2H}\,.
\end{equation}
In addition, we can take $\pier{\rho}= -\Delta (\psi_n(t) ) $ in \eqref{ls2n} and integrate by parts. We then obtain
\begin{align}
& \|\Delta \psi_n(t)\|^2_H \,=\, 
(\tau \dt \psi_n (t)+ f''(\phis(t))\psi_n(t) - v_h(t),
\Delta \psi_n(t)) + \iO \nabla \eta_n (t) \cdot \nabla \psi_n(t) \non\\
& \leq \frac 1 2 \|\Delta \psi_n(t)\|^2_H  + 
C \bigl(1+ \|\psi_n\|_{\W{1,\infty}H}^2  +  \|v_h\|_{\L{\infty}H}^2 \bigr)
\non \\ 
&\quad{}+ \|\nabla \eta_n \|_{{\L{\infty}H}^3} \|\nabla \psi_n \|_{{\L{\infty}H}^3}
\quad \hbox{for a.e. } t\in (0,T)\,.
\label{estP3}
\end{align}
Consequently, from \eqref{esti2} and the elliptic regularity theory,
we find that 
\begin{equation}
\label{estP4}
 \| \Delta \psi_n \|_{\L{\infty}H} +  \| \psi_n \|_{\L{\infty}W} \,\le\,C\|h\|_{\L2H}\,.
\end{equation}

\step
Conclusion of the proof

By virtue of the estimates shown above, we conclude that, for all $n\in\enne$,
\begin{align}
\label{lsboundn}
&\|\eta_n\|_{W^{2,\infty}(0,T;\VD)\cap W^{1,\infty}(0,T;H)\cap L^\infty(0,T;V)}
+ \|\psi_n\|_{W^{1,\infty}(0,T;H)\cap H^1(0,T;V)\cap L^\infty(0,T;W)}\non\\
&\le\,C\|h\|_{\L2H}\,. 
\end{align} 
Hence there exists a pair $(\eta,\psi)$ such that (possibly on a subsequence, 
which is again labeled by $n\in\enne$)
\begin{align}
\label{coneta}
&\eta_n\to\eta \quad \mbox{weakly star in }\,\W{2,\infty}{V^*}\cap \W{1,\infty}H \cap \L\infty V,\\
\label{conpsi}
&\psi_n\to\psi \quad \mbox{weakly star in } \,W^{1,\infty}(0,T;H)\cap H^1(0,T;V)\cap L^\infty(0,T;W).
\end{align}
By a standard argument (which needs no repetition here) it then follows that the pair $(\eta,\psi)$ solves the equations
\eqref{ls1}, \eqref{ls2} with $v=v_h$, and \eqref{ls4}. In addition, we easily conclude that $\eta(0)=\psi(0)=0$. Moreover, it follows from
\eqref{lsboundn} and the lower semicontinuity properties of norms that the triple $(\eta_h,\psi_h,v_h)$, where $\eta_h=\eta$ and 
$\psi_h=\psi$\pier{,} satisfies the estimate \eqref{contin}. From this it is easily seen that the solution is unique: indeed, if $(\eta_i,\psi_i,v_i)$,
$i=1,2$, are two solutions, then the triple $(\eta,\psi,v):=(\eta_1-\eta_2,\psi_1-\psi_2,v_1-v_2)$ solves the system \LinState\ with
$h=0$, whence $(\eta,\psi,v)=(0,0,0)$ follows. With this, the assertion is completely proved. 
\Edim

Next, we show a differentiability result for the control-to-state operators $\CS^\alpha$, for $\alpha\in (0,1]$.
\Bthm
\label{tre.due}
Suppose that {\bf (A1)}--{\bf (A5)} are fulfilled, and let $\alpha\in (0,1]$ be given. Then the control-to-state operator $\CS^\alpha$ 
is Fr\'echet differentiable in $\CU_R$ as a mapping from $\LiQ$ into the Banach space 
\begin{equation}
\label{defX}
{\cal X}:= C^0([0,T];H)\times \bigl(\H1H\cap C^0([0,T];V)\cap \L2W\bigr)\times\H1H.
\end{equation}
Moreover, for every $u\in\CU_R$ the Fr\'echet derivative $D\CS^\alpha(u)\in {\cal L}(\LiQ,{\cal X})$ is defined as follows:
for all $h\in\LiQ$ it holds \,$D\CS^\alpha(u)[h]=(\eta_h,\psi_h,v_h)$, where $(\eta_h,\psi_h,v_h)$ is the unique solution to the 
linearized system \LinState.   
\Ethm
\Bdim
Let $\alpha\in(0,1]$ and $u\in\CU_R$ be given, and let $(\mus,\phis,\ws)=\CS^\alpha(u)$. 
Since $\CU_R$ is open, it holds
$u+h\in\CU_R$ provided $\|h\|_{\LiQ}$ is sufficiently small. In the following we only consider such variations $h$, and we define
the quantities 
\begin{align*}
&(\mu_h,\phi_h,w_h):=\CS^\alpha(u+h), \quad y_h:=\mu_h-\mus-\eta_h,\non\\
& z_h:=\phi_h-\phis-\psi_h, \quad \omega_h:=w_h-\ws-v_h,
\end{align*}
where $(\eta_h,\psi_h,v_h)$ denotes the unique solution to the linearized system \LinState. Notice that 
\begin{align*}
&y_h\in W^{2,\infty}(0,T;\VD)\cap W^{1,\infty}(0,T;H)\cap \cap L^\infty(0,T;V), \non\\
&z_h \in W^{1,\infty}(0,T;H)\cap H^1(0,T;V)\cap L^\infty(0,T;W), \non\\
&\omega_h\in H^1(0,T;H),
\end{align*}
according to Theorem~\ref{due.uno} and Theorem~\ref{tre.uno}. Moreover, from the estimate \eqref{stimaK4} in Theorem~\pier{\ref{due.tre}} we conclude that
\begin{align}
&\norma{\mu_h-\mus}_{\H2\VD\cap\W{1,\infty} H \cap \L\infty V}
+ \norma{\phi_h-\phis}_{\H1H\cap\L\infty V\L2W}
\non\\
&\le \,C \norma{h}_{\L2H}\,,
  \label{elvis1}
\end{align}  
where, here and in the remainder of this proof, $C>0$ denotes constants that may depend on $\alpha, \,\,R$, and the data of the
system, but not on the special choice of $h$ with $u+h\in\CU_R$.  

Next, we observe that it follows from the estimate \eqref{contin} in Theorem~\ref{tre.uno} 
that the mapping $\,h\mapsto (\eta_h,\psi_h,v_h)\,$
belongs to ${\cal L}(\LiQ,{\cal X})$, in particular. Hence it suffices to show that
\begin{equation}
\label{elvis2}
\lim_{\|h\|_{\LiQ}\to 0} \frac{\|(y_h,z_h,\omega_h)\|_{\cal X}}{\|h\|_{\LiQ}}=0,
\end{equation}
which is certainly satisfied if 
\begin{align}
\label{reicht}
&\|y_h\|^2_{C^0([0,T];H)}\,+\,\|z_h\|^2_{\H1H \cap C^0([0,T];V)\cap \L2W}\,+\,\|\omega_h\|^2_{\H1H}\,\le\,C\|h\|_{\L2H}^4 
\end{align}
for all $h\in\LiQ$ with $u+h\in\CU_R$. We are going to show that \eqref{reicht} is in fact valid. To this end, observe that
the triple $(y_h,z_h,\omega_h)$ apparently solves the system\pier{%
\begin{align}
\label{diffi1}
&\alpha\langle\dtt y_h,\rho\rangle + ( \dt z_h,\rho)+\iO\nabla y_h \cdot \nabla \rho   = 0\quad \hbox{for every  $\rho\in V$}, \, \hbox{ a.e. in $(0,T)$},\\
\label{diffi2}
&  \tau \dt z_h - \Delta z_h = y_h+\omega_h - \bigl(f'(\phi_h)-f'(\phis)-f''(\phis)\psi_h\bigr)\quad \hbox{a.e. in $Q$},\\
\label{diffi3}
&\gamma\dt\omega_h+\omega_h=0\quad \mbox{a.e. in $Q$,}\\
\label{diffi4}
&\dn z_h=0\quad \mbox{a.e. on $\Sigma$},\\
\label{diffi5}
&y_h(0)=\dt y_h(0)=z_h(0)=\omega_h(0)=0 \quad \mbox{a.e. in $\Omega$.}
\end{align}
Obviously, it results that}  $\omega_h=0$ a.e. in $Q$. Moreover, Taylor's theorem with integral remainder yields almost everywhere in $Q$  the identity
\begin{align}
\label{Taylor}
&f'(\phi_h)-f'(\phis)-f''(\phis)\psi_h=f''(\phis)z_h + A_h\,(\phi_h-\phis)^2,\quad\mbox{where}\non\\
&A_h:=\int_0^1\int_0^1 s\,f'''(\phis+\tau\,s\,(\phi_h-\phis))\,d\tau\,ds\,.
\end{align}
Owing to {\bf (A3)} and \eqref{separa}, we have
\begin{equation}
\label{estiAh}
|A_h|\,\le\,C \quad\mbox{a.e. in } Q.
\end{equation}
We now integrate \eqref{diffi1} over $(0,t)$ for $t\in(0,T]$ to obtain the equation
\begin{equation}
\label{diffi1neu}
\alpha(\dt y_h,\rho)+(z_h,\rho)+\iO \nabla(1\ast y_h)\cdot\nabla\rho=0,
\end{equation}
where the expression $1\ast y_h$ is defined in \eqref{anti}. We test \eqref{diffi1neu} by $y_h$, 
\eqref{diffi2} by $z_h$, add the resulting identities, and integrate with respect to time over $(0,t)$ where 
$t\in(0,T]$. Noting the cancellation of two terms, we then arrive at the identity
\begin{align}
\label{elvis3}
&\frac \alpha 2 \|y_h(t)\|_H^2 \,+\,\frac 12 \|\nabla(1\ast y_h(t))\|^2_{H\times H\times H}\,+\,\frac \tau 2\|z_h(t)\|_H^2
\,+\iint_{Q_t} |\nabla z_h|^2 \non\\
&= -\iint_{Q_t} f''(\phis)|z_h|^2 \,-\iint_{Q_t}A_h(\phi_h-\phis)^2z_h\,.
\end{align}
Thanks to {\bf (A3)}, \eqref{separa}, \eqref{elvis1},  and \eqref{estiAh}, the \rhs\ of this equation can be bounded by
\begin{align}
\label{elvis4}
&C\iint_{Q_t}|z_h|^2\,+\,C\iot\|z_h(s)\|_H\,\|\phi_h(s)-\phis(s)\|_{L^4(\Omega)}^2\,ds\non\\
&\le\,C\iint_{Q_t}|z_h|^2\,+\,C\iot \|\phi_h(s)-\phis(s)\|^4_{L^4(\Omega)}\,ds\non\\
&\le\,C\iint_{Q_t}|z_h|^2\,+\,C\|h\|_{\L2H}^4\,,
\end{align}
where in the last estimate the continuity of the embedding $V\subset L^4(\Omega)$ was used. Combining \eqref{elvis3} with   
\eqref{elvis4}, and invoking Gronwall's lemma, we then conclude that 
\begin{align}
\label{reicht1}
&\|y_h\|^2_{C^0([0,T];H)}\,+\,\|z_h\|^2_{C^0([0,T];H)\cap \L2V}\,\le\,C\|h\|^4_{\L2H}\,.
\end{align}
 
At this point we note that $z_h$ solves a linear parabolic initial-boundary value problem of heat conduction type, 
with zero initial datum and zero Neumann boundary 
condition, whose right-hand side $g_h:=y_h-f''(\phis)z_h -A_h(\phi_h-\phis)^2$ has already been shown to satisfy the condition
\,$\|g_h\|_{\L2H}^2\,\le\,C\|h\|_{\L2H}^4$. It therefore follows from standard linear parabolic theory that
\begin{align}
\label{reicht2}
\|z_h\|^2_{\H1H \cap C^0([0,T];V)\cap \L2W}\,\le\,C\|h\|^4_{\L2H}\,.
\end{align}
Combining \eqref{reicht1} with \eqref{reicht2}, and recalling that $\omega_h=0$, we obtain that \eqref{reicht} is valid, 
which finishes the proof of the assertion.
\Edim   

\Brem
Note that the value of the constant $R>0$ did not really matter in the above proof. We therefore conclude that the 
Fr\'echet derivatives $D\CS^\alpha(u)\in{\cal L}(\LiQ,{\cal X})$ exist for all $u\in\LiQ$.
\Erem

\section{The optimal control problem}
\setcounter{equation}{0}

In this section, we investigate the control problem formulated in the introduction, which, in order to stress its dependence
on the parameter $\alpha\in(0,1]$, will in the following be denoted by {\bf (CP$_\alpha$)}. In addition to the general 
assumptions {\bf (A1)}--{\bf (A5)},
we postulate:

\vspace*{2mm} \noindent
{\bf (A6)}  \quad  $b_1\ge 0$, $b_2\ge 0$, $b_3>0$ and $\kappa\ge 0$ are given constants.

\vspace*{1mm} \noindent
{\bf (A7)} \quad $\phi_Q\in L^2(Q)$, $\phi_\Omega\in V$, and $\underline u,\overline u\in\LiQ$ are given functions 
such that $\,\underline u\le\overline u\,$ almost everywhere in $Q$.

\vspace*{1mm} \noindent
{\bf (A8)} \quad The mapping $G:L^2(Q)\to\erre$ is convex, continuous, and nonnegative. 

\subsection{Existence}
We begin our analysis with an existence result.
\Bthm
\label{quattro.uno}
Suppose that {\bf (A1)}--{\bf (A8)} are fulfilled. Then the optimal control problem {\bf (CP$_\alpha$)} has for every 
$\alpha\in (0,1]$ at least one solution.
\Ethm
\Bdim Obviously, the cost functional $\CJ$ is nonnegative. Therefore, there is a minimizing sequence $\{u_n\}_{n\in\enne}
\subset \Uad$, i.e., we have
$$
\lim_{n\to\infty}\CJ(\CS_2^\alpha(u_n),u_n)=\inf_{u\in\Uad}\CJ(\CS^\alpha_2(u),u)\ge 0.
$$ 
Then the triple $(\mu^{\alpha}_n,\phi^\alpha_n,w^\alpha_n):=\CS^\alpha(u_n)$ solves for $n\in\enne$ the state system \State\ with $u=u_n$ and thus
satisfies the global bounds \eqref{stimaK1} and \eqref{ssboundal}. We may therefore assume without loss of generality that,
with a suitable triple $(\mu,\phi,w)$,
\begin{align}
\label{conmun}
\mu^\alpha_n\to\mu &\quad\mbox{weakly star in } W^{2,\infty}(0,T;V^*) \cap W^{1,\infty}(0,T;H) \cap L^\infty(0,T; V)\,,\\
\label{conphin}
\phi^\alpha_n\to\phi&\quad \mbox{weakly star in } \W{1,\infty}H \cap \H1V \cap \L\infty W\,,\\
\label{conwn}
w^\alpha_n\to w&\quad \mbox{weakly star in } W^{1,\infty}(0,T;\pier{\Linfty})\,,\\
\label{gleichm}
\phi^\alpha_n\to\phi&\quad\mbox{strongly in } C^0(\overline Q)\,,
\end{align}
where the last convergence result is a consequence of \cite[Sect.~8, Cor.~4]{Simon} and the compactness of the embedding 
$W\subset C^0(\overline\Omega)$. In addition, we have $\,r_*(\alpha)\le\phi^\alpha_n\le r^*(\alpha)\,$ in $\overline Q$ for every $n\in\enne$.
Hence, \eqref{gleichm} and {\bf (A3)} imply that 
\begin{equation} \label{confn}
f'(\phi^\alpha_n)\to f'(\phi) \quad\mbox{strongly in } C^0(\overline Q).
\end{equation}
\pier{Moreover,}, $\Uad$ is a convex, bounded and closed subset of $L^\infty(Q)$. Therefore, we may assume that there is some $u^\alpha\in\Uad$ such that
\begin{equation}
\label{conun}
u_n\to u^\alpha\quad\mbox{weakly star in }\LiQ.
\end{equation}
Passage to the limit as $n\to\infty$ in the equations \State, written for $((\mu^\alpha_n,\phi^\alpha_n, w^\alpha_n), u_n)$, taking 
\eqref{conmun}--\eqref{conun} 
into account, then leads to the conclusion that $(\mu,\phi,w)$ solves the system \State\ for the control $u=\us$. In other words,
we have $(\mu,\phi,w)=(\mus,\phis,\ws)=\CS^\alpha(\us)$. It then follows from the semicontinuity properties of the functionals $J$ and $G$ 
(notice that $G$ is convex and continuous and thus weakly sequentially lower semicontinuous on $L^2(Q)$) that
$$\CJ(\phis,\us)=J(\phis,\us)+\kappa G(\us)\le \liminf_{n\to\infty}(J(\phi^\alpha_n,u_n)+\kappa G(u_n))=\inf_{u\in\Uad} \CJ(\CS_2^\alpha(u),u),$$
which means that $u=\us$ is an optimal control.
\Edim 

\subsection{Necessary optimality conditions} 
We begin by deriving a first condition for a control to be locally optimal. To this end, recall that $\us\in\Uad$ is called locally
optimal for \CP\ in the sense of $L^p(Q)$ for some $p\in[1,+\infty]$ if and only if there is some $\varepsilon>0$ such that
$\CJ(\CS^\alpha_2(\us),\us)\,\le\,\CJ(\CS^\alpha_2(u),u)$ for all $u\in\Uad$ with \,$\|u-\us\|_{L^p(Q)}\le \varepsilon$.
Observe that every control which is locally optimal in the sense of $L^p(Q)$ for some $p\in[1,+\infty)$ is also locally optimal in the sense of 
$\LiQ$. Therefore, all necessary conditions derived for locally optimal controls in the sense of $\LiQ$ are automatically valid also for
locally optimal controls in the sense of $L^p(Q)$ for any $p\in[1,+\infty)$ and thus, in particular, for (globally) optimal controls.

We  have the following result.
\Blem
Suppose that {\bf (A1)}--{\bf (A8)} are fulfilled and $\alpha\in(0,1]$. Moreover, assume that $\us\in\Uad$ is a locally optimal
control for \CP\ in the sense of $\LiQ$ with associated state $(\mus,\phis,\ws)=\CS^\alpha(\us)$. Then there is some $\lambda^\alpha\in\partial G(\us)$ 
(where $\partial G(\us)\subset L^2(Q)$ denotes the subdifferential of $G$ at $\us$) such that it holds
\begin{align}
\label{VUG1}
&b_1\iint_Q 	\psi_h (\phis-\phi_Q)\,+\,b_2\iO \psi_h(T)(\phis(T)-\phi_\Omega)\,+\,b_3\iint_Q \us(u-\us)\non\\
&+\,\kappa\iint_Q\lambda^\alpha(u-\us)\,\ge\,0 \quad\mbox{for all } u\in\Uad,
\end{align}
where $(\eta_h,\psi_h,v_h)$ is the unique solution to the linearized system \LinState\ associated with $h=u-\us$.
\Elem
\Bdim
For the following, we introduce the reduced functionals $\widehat J^\alpha$ and $\widehat \CJ^\alpha$ by setting
\begin{equation}
\label{reduced}
\widehat J^\alpha(u):=J(\CS_2^\alpha(u),u), \quad \widehat \CJ^\alpha(u):=\widehat J^\alpha(u)+\kappa G(u), \quad\mbox{for }u\in\LiQ.
\end{equation}
Then, by Theorem~\ref{tre.due} and the quadratic form of the differentiable part $J$ of the cost functional, it follows from the chain rule
 that the mapping $\widehat J^\alpha:\LiQ\to \erre$ has a Fr\'echet derivative $D\widehat J^\alpha(u)$ for every $u\in\LiQ$. From the
convexity of $\Uad$ and of $G$ we then infer that 
$$D\widehat J^\alpha(\us)[u-\us]+\kappa(G(u)-G(\us))\,\ge\,0\quad\mbox{for all } u\in\Uad.$$
A standard argument from convex analysis (for the details, see, e.g., \cite[Sect.~4.1]{STAllenCahn}) yields that there is some $\lambda^\alpha
\in \partial G(\us)$ such that 
\begin{equation}
\label{uffuff}
D\widehat J^\alpha(\us)[u-\us]+\kappa\iint_Q\lambda^\alpha(u-\us)\,\ge\,0 \quad\mbox{for all } u\in\Uad.
\end{equation}
Finally, we infer from the chain rule, using Theorem~\ref{tre.due} and the special form of $J$ (which, in particular, does not involve the 
first and third solution components), that 
\begin{align*}
&D\widehat J^\alpha(\us)[u-\us]\\
&{}=b_1\iint_Q \!	\psi_h (\phis-\phi_Q)\,+\,b_2\iO \!\psi_h(T)(\phis(T)-\phi_\Omega)\,+\,b_3\iint_Q \!\us(u-\us),
\end{align*}
where $h=u-\us$. This concludes the proof of the assertion.
\Edim

As usual, we improve the still rather useless optimality condition \eqref{VUG1} by means of the adjoint variables. 
We write the associated adjoint system
in the form
\begin{align}
\label{as2}
&\pier{{}\alpha \dtt p-\Delta p - q=0}&&\mbox{a.e. in }Q\,,\\
\label{as1}
&\pier{{}-\dt(p+\tau q)-\Delta q+f''(\phis)q=b_1(\phis-\phi_Q)}&&\mbox{a.e. in }Q\,,\\
\label{as3}
&-\gamma\dt r+r-q=0 &&\mbox{a.e. in }Q\,,\\
\label{as4}
&\dn p=\dn q =0 &&\mbox{a.e. on }\Sigma\,,\\
\label{as5}
&p(T)=\dt p(T)=r(T)=0, \,\,\, q(T)=\frac{b_2}{\tau}(\phis(T)-\phi_\Omega)&&\mbox{a.e. in }\Omega\,.
\end{align}

We have the following well-posedness result.
\Bthm
\label{quattro.tre}
Suppose that {\bf (A1)}--{\bf (A8)} are fulfilled and that $\alpha\in (0,1]$. In addition, let $\us\in\CU_R$ be given and 
$(\mus,\phis,\ws)=\CS^\alpha(\us)$. Then the system \pier{\eqref{as2}--\eqref{as5}} has a unique solution $(\ps,\qs,\rs)$ with the
regularity
\begin{align}
\label{regp}
&\ps\in W^{2,\infty}(0,T;H)\cap W^{1,\infty}(0,T;V)\cap \L\infty W\,,\\
\label{regq}
&\qs\in \H1H \cap C^0([0,T];V)\cap \L2 W\,,\\
\label{regr}
&\rs\in \H2H \cap C^1([0,T];V) \cap \H1 W\,.
\end{align}
Moreover, there is a constant \pier{$K_5(\alpha) >0$}, which depends only on $R$ and the data of the state system, such that
\begin{align}
\label{stimaAdj1}
&\alpha\|\ps\|_{W^{2,\infty}(0,T;H)}+\alpha^{1/2} \|\ps\|_{W^{1,\infty}(0,T;V)}+\|\ps\|_{\L\infty W}
\non\\
&\pier{{}+\|\qs\|_{\L\infty V \cap\L2W}+\|\rs\|_{W^{1,\infty}(0,T;V)\cap \H1 W}}\non \\
&\pier{{}+\|\ps+\tau\qs\|_{\H1 H} +\alpha^{1/2} \|\dt \qs\|_{\L2H}+ \alpha^{1/2} \|\dtt \rs\|_{\L2 H}\,\le\,K_5(\alpha).}
\end{align}
\Ethm
\Bdim
The existence proof follows that of Theorem~\ref{tre.uno}: we again employ a Faedo--Galerkin approximation using the same basis functions as
in the proof of Theorem~\ref{tre.uno}. To keep the paper at a reasonable length, we avoid here to write the approximating
system explicitly and only give formal estimates for $(p,q,r)$ that would be rigorous on the level of the finite-dimensional
approximations. In the following, we use the notation
$$Q^t:=\Omega\times (t,T)\pier{,}\quad\mbox{for } t\in [0,T).$$

\step First estimate

To begin with, we first test \pier{\eqref{as2} by $\,-\dt p$, \eqref{as1} by $q$,} and \eqref{as3} by $\,-\dt r$, and add
the resulting equations, noting the cancellation of two terms. Then we integrate over $(t,T)$, where $t\in[0,T)$ is arbitrary.
Using the terminal conditions \eqref{as5}, we arrive at the identity  
\begin{align}
&\frac \alpha 2\|\dt p(t)\|_H^2 + \frac 12 \iO|\nabla p(t)|^2 + \frac \tau 2 \|q(t)\|_H^2 +\iint_{Q^t}|\nabla q|^2
+\gamma\iint_{Q^t}|\dt r|^2 +  \frac 12 \|r(t)\|^2_H \non \\
&= \frac{b_2^2}{2\tau}\iO|\phis(T)-\phi_\Omega|^2-\iint_{Q^t} f''(\phis)q^2+\iint_{Q^t}\bigl(b_1(\phis-\phi_Q)-\dt r\bigr)q\,. \label{pier1} 
\end{align}
\pier{The first term on the right-hand side of \eqref{pier1} is under control due to \eqref{stimaK1} and {\bf (A7)}. Next, we note that
$$-\iint_{Q^t} f''(\phis)q^2 = -\iint_{Q^t} f_1''(\phis)q^2 -\iint_{Q^t} f_2''(\phis)q^2, $$
with the first contribution being nonpositive and the second one being bounded by $C \iint_{Q^t} q^2$, due to~{\bf (A2)}, \eqref{stimaK1} and {\bf (A3)}. About the third term on the right-hand side of \eqref{pier1}, by \eqref{stimaK1}, {\bf (A7)} and Young's inequality we have that
$$ \iint_{Q^t}\bigl(b_1(\phis-\phi_Q)-\dt r\bigr)q \leq \frac \gamma 2 \iint_{Q^t}|\dt r|^2 + C \iint_{Q^t} q^2.$$
Then,} using Gronwall's lemma, we immediately see that
\begin{align}
\label{estad1}
&\alpha^{1/2}\|\dt p\|_{\L\infty H}+\|\nabla p\|_{{\L\infty H}^3}+\|q\|_{\L\infty H \cap \L2V} +\|r\|_{\H1H}\,\le\,C\,,
\end{align}
where \pier{the constant $C>0$ depends on the data of the system, but is independent of}~$\alpha\in(0,1]$.

\step Second estimate

Next, we test \eqref{as1} by $p+\tau q$, obtaining the identity
\begin{align*}
&\frac 12 \|(p+\tau q)(t)\|_H^2 +\iint_{Q^t}\nabla q\cdot\nabla(p+\tau  q)\\
&=\frac{\tau^2}2\|q(T)\|^2 +\iint_{Q^t} \bigl(b_1(\phis-\phi_Q)-f''(\phis)q\bigr)(p+\tau q)\,.
\end{align*}
Hence, by virtue of Young's inequality\pier{, \eqref{estad1} and \eqref{ssboundal}}, we find from Gronwall's lemma that $\,\|p+\tau q\|_{\L\infty H}\,\le\,\pier{C_\alpha}$. \pier{Therefore, using once more \eqref{estad1}, it turns out that \,$\|p\|_{\L\infty H}\,\le\,C_\alpha$,} which, in turn, implies that
\begin{equation}
\label{estad2}
\|p\|_{\L\infty V}\,\le\,\pier{C_\alpha}\,.
\end{equation}

\step Third estimate

In this estimate, which is entirely formal (but completely justified on the level of the Faedo--Galerkin approximations), we first test 
\eqref{as2} by $\Delta\dt p$. Taking the boundary conditions \eqref{as4} and the terminal conditions \eqref{as5} into account,
we obtain the identity
\begin{equation}\label{Kurt1}
\frac{\alpha}2 \|\nabla\dt p(t)\|_H^2 +\frac 12 \|\Delta p(t)\|_H^2=\iint_{Q^t}q \,\Delta\dt p=\iint_{Q^t}\Delta q\,\dt p\,.
\end{equation}
\pier{Secondly}, we test \eqref{as1} by \,$-\Delta q$\, to see that
\begin{align}\label{Kurt2}
&\iint_{Q^t} \Delta q \,\dt p + \frac\tau 2\iO|\nabla q(t)|^2+\iint_{Q^t}|\Delta q|^2\non\\
&= \frac \tau   2 \iO |\nabla q(T)|^2 + \iint_{Q^t}\Delta q\bigl(f''(\phis)q-b_1(\phis-\phi_Q)\bigr)\,.                 
\end{align}
\pier{In view of \eqref{as5}}, we recall that $\,\phi_\Omega\in V$, by assumption {\bf (A7)} \pier{and consequently, by \eqref{stimaK1}, deduce that}
$$\iO|\nabla q(T)|^2\,\le\,C(\|\phis(T)\|_V^2+\|\phi_\Omega\|_V^2)\le C\,.$$
Thus, adding \eqref{Kurt1} and \eqref{Kurt2} \pier{leads to a cancellation of two terms. Using \eqref{ssboundal} and Young's inequality, by} standard elliptic estimates we can infer that
\begin{align}
&\alpha^{1/2}\norma{p} _{\W{1,\infty}V} + \norma{p}_{\L\infty W}+\|q\|_{\L\infty V\cap \L2W}\,\le\,\pier{C_\alpha}\,.
\label{estad3}
\end{align}
\pier{It} then directly follows from \eqref{as2} and \eqref{estad3} that
$
\|\alpha\dtt p\|_{\L\infty H}\,\le\,\pier{C_\alpha}\,,
$ 
and thus, by virtue of \eqref{estad1},
\begin{equation}
\label{estad4}
\pier{\alpha\| p\|_{W^{2,\infty}(0,T;H)}\,\le\,C_\alpha}\,.
\end{equation}
In addition, \pier{a comparison of terms in \eqref{as1} yields}
\begin{equation}
\label{estad5}
\|p+\tau q\|_{\H1H}\,\le\,\pier{C_\alpha}\,.
\end{equation}
\pier{Moreover, due to the estimate for $q$ in \eqref{estad3},  arguing on \eqref{as3} and using the terminal condition for $r$ in \eqref{as4} enable us to conclude that}
\begin{equation}
\label{estad6}
\|r\|_{W^{1,\infty}(0,T;V)\cap \H1W}\,\le\,\pier{C_\alpha}\,.
\end{equation}
\pier{The last part regards a similar bound for 
$$\alpha^{1/2} \|\dt \qs\|_{\L2H}+ \alpha^{1/2} \|\dtt \rs\|_{\L2 H}$$
which follows then from \eqref{estad1}, \eqref{estad5} and comparison in \eqref{as1} and \eqref{as3}. 
At this point, we have completely shown the relevant estimate \eqref{stimaAdj1} that suffices to prove} the existence of a solution to the adjoint system
\pier{\eqref{as2}--\eqref{as5}} that has the asserted properties. Owing to the linearity of the system, the proof of uniqueness is 
simple and can be skipped here. 
\Edim

\Brem
\label{PIER2}
\pier{Recalling the contents of Remark~\ref{PIER1} and the assumption~{\bf (A3)}, we aim to point out that 
in the case when  $\,D(\partial f_1)= \erre$\, and \,$f_1 \in C^3 (\erre)$\, in {\bf (A3)}, then the constant $K_5(\alpha)$ in \eqref{stimaAdj1} may be chosen independent of $\alpha \in (0,1]$,  since in this framework the coefficient $f''(\phis)$ appearing in \eqref{as1} is bounded in $L^\infty (Q) $ independently of $\alpha$. In fact, we can repeat all the argumentation leading to \eqref{estad2}--\eqref{estad6} with constants $C$ in the \rhs s that do not depend on $\alpha$.}
\Erem

Having shown the well-posedness of the adjoint system, it is now a standard procedure to improve the necessary condition 
for locally optimal controls stated in Lemma 4.2. Indeed, a straightforward calculation using the linearized system \LinState\
and the adjoint system \pier{\eqref{as2}--\eqref{as5}}, which can be left to the reader, leads to the
following result.
\Bthm
\label{quattro.cinque}
Suppose that {\bf (A1)}--{\bf (A8)} are fulfilled and $\alpha\in(0,1]$. Moreover, assume that $\us\in\Uad$ is a locally optimal
control for {\bf (CP$_\alpha$)} in the sense of $\LiQ$ with associated state $(\mus,\phis,\ws)=\CS^\alpha(\us)$ and adjoint state $(\ps,\qs,\rs)$.
Then there is some $\lambda^\alpha\in\partial G(\us)$ such that
\begin{align}
\label{VUG2}
&\iint_Q 	(\rs\,+\,b_3\us \,+\,\kappa\lambda^\alpha)(u-\us)\,\ge\,0 \quad\mbox{for all } u\in\Uad\,.
\end{align}
\Ethm

\subsection{Sparsity of controls}

The convex functional $G$ in the cost functional accounts for the sparsity of optimal controls, i.e., the possibility that
every locally optimal control may vanish in some subset of the space-time cylinder $Q$. The form of this region depends on the choice 
of $G$, where the sparsity properties can be deduced from the variational inequality \eqref{VUG2} and the form of the 
subdifferential $\partial G$. In the following, we restrict ourselves to the case of {\em full sparsity} which is obtained for the  
case of the $L^1(Q)$-norm \eqref{defG} that recently has been investigated in \cite{CoSpTr}. Arguing exactly as there,
one obtains the following result (see~\cite[Thm.~4.7]{CoSpTr}).
\Bthm
\label{quattro.sei}
Suppose that {\bf (A1)}--{\bf (A8)} are fulfilled, $\alpha\in (0,1]$, and $\kappa>0$. Assume that $\underline u$ and $\overline u$ are constants such that
$\underline u<0<\overline u$. If $\us\in\Uad$ is locally optimal in the sense of $\LiQ$ for {\bf (CP$_\alpha$)} with associated 
state $(\mus,\phis,\ws)$ and adjoint state $(\ps,\qs,\rs)$, then there exists a function $\lambda^\alpha\in\partial G(\us)$ satisfying \eqref{VUG2},
and it holds

\begin{equation}
\us(x,t)=0 \quad\mbox{if and only if}\quad |\rs(x,t)|\le\kappa, \quad\mbox{for a.e. }(x,t)\in Q\,.
\end{equation}
Moreover, if $\rs$ and $\lambda^\alpha$ are given, then $\us$ is obtained from the projection formula
$$\us(x,t)\,=\,\max \left \{  \underline u,\,\min \left \{\overline u,\,-b_3^{-1}(\rs+\kappa\lambda^\alpha)(x,t) \right\} \right\} 
\quad\mbox{for a.e. $(x,t)\in Q$}.$$
\Ethm

\Brem
Owing to the global estimate \eqref{stimaAdj1}, we have 
$$
%\sup_{\alpha\in (0,1]} 
\pier{\|\rs\|_{\H1 W}\,\le\,\pier{K_5(\alpha)}},$$
and thus, thanks to the continuity of the embedding $\H1W \subset C^0(\overline Q)$, 
$$
%\sup_{\alpha\in (0,1]}
\pier{\|\rs\|_{C^0(\overline Q) }\,\le\,\widehat \kappa (\alpha)},$$
with a constant $\pier{\widehat\kappa (\alpha)}>0$. \pier{Hence, for each  $\alpha\in (0,1]$, if $\kappa\ge\widehat \kappa (\alpha)$ then every} 
locally optimal control in the sense of $\LiQ$ for {\bf (CP$_\alpha$)} must vanish.  
\Erem

\section{Asymptotic analysis}
\label{Asy}
\setcounter{equation}{0}

This section is devoted to the study of the asymptotic behavior of the problem {\bf (CP$_\alpha$)} as $\alpha\searrow0$. 
The case $\alpha=0$, which has
been thoroughly investigated in \cite{CoSpTr}, is denoted by {\bf (CP$_0$)}. 

\subsection{Convergence of the state variables}
The asymptotic results for the state system have already been studied in \cite{CS12} where also the case of the nondifferentiable double 
obstacle potential \eqref{obspot} was included in the analysis (see, \cite[Thm.~5.1]{CS12}). We have the following result.
\Bthm
\label{connystate} 
Assume that {\bf (A1)}--{\bf (A5)} are satisfied, and let $\{\alpha_n\}_{n\in\enne}\subset (0,1]$ and $\{u^{\alpha_n}\}_{n\in\enne}
\subset\Uad$
be sequences such that $\,\alpha_n\searrow 0$\, and $\,u^{\alpha_n}\to u^0\,$ weakly star in $\LiQ$.  
If  $(\mu^{\alpha_n} , \phi^{\alpha_n} , w^{\alpha_n}):=\CS^{\alpha_n}(u^{\alpha_n})$ denotes for $n\in\enne$ the unique 
solution to the state system \State\ for $u=u^{\alpha_n}$ established in Theorem~\ref{due.uno}, then it holds
\begin{align}
\label{connymu1}
&\mu^{\alpha_n} \to\mu^0 \quad \mbox{weakly star in }\,\pier{\L\infty V},\\
\label{connymu2}
&\alpha_n\,\mu^{\alpha_n} \to 0 \quad \mbox{weakly star in }\,\pier{\W{2,\infty}{V^*}}\,  
\mbox{ and strongly in }\,\pier{\W{1,\infty}H},\\
\label{connyphi}
&\phi^{\alpha_n}\to\phi^0 \quad \mbox{weakly star in } \,W^{1,\infty}(0,T;H)\cap H^1(0,T;V)\cap L^\infty(0,T;W)
\non
\\
&\qquad\qquad \quad \mbox{ and strongly in }\,\C0V \cap C^0(\overline Q),\\
&\pier{f'(\phi^{\alpha_n}) \to f'(\phi^0) \quad \mbox{weakly star in }  L^\infty(0,T;H) \mbox{ and a.e.~in }\,Q,}
\label{confpphi}
\\
&w^{\alpha_n}\to w^0 \quad \mbox{weakly star in } W^{1,\infty}(0,T;\pier{\Linfty}), \label{connyw}
\end{align}
where $(\mu^0,\phi^0,w^0)$ is the unique strong solution to the viscous Cahn--Hilliard system 
\begin{align}
\label{ss01}
& \dt\phi-\Delta\mu  = 0 &&\mbox{a.e. in }Q\,,\\  
\label{ss02}
&  \tau \dt \phi- \Delta \phi   + f'(\phi) = \mu + w &&\mbox{a.e. in }Q\,,\\ 
\label{ss03}
&\gamma \dt w+w=u^0&&\mbox{a.e. in }Q\,,\\
\label{ss04}
& \dn \mu=\dn\phi=0 && \mbox{a.e. on }\Sigma\,,\\
\label{ss05}
& \mu(0)=\mu_0, \quad \phi(0)= \phi_0,\quad w(0)=w_0 && \hbox{a.e.~in }{\Omega} \,.
\end{align}
\Ethm

\Bdim
\pier{Arguing as in the proof of Theorem~\ref{quattro.uno}, by \eqref{stimaK1} it is not difficult to check}
that $(\mu^{\alpha_n},\phi^{\alpha_n},w^{\alpha_n})$ converges to a triple $(\mu,\phi,w)$ in the sense of 
\pier{\eqref{connymu1}--\eqref{connyphi}, \eqref{connyw},} at least on a subsequence which we still index by $n\in\enne$. In view of \eqref{connyphi} \pier{and the Lipschitz continuity of $f_2'$,  we have that $f_2'(\phi^{\alpha_n})\to f_2'(\phi)$ strongly in $\,C^0(\overline Q)$. Moreover, recalling \eqref{stimaK1}, {\bf (A3)} and 
the estimate~\eqref{ssboundal}, we also have that $f_1'(\phi^{\alpha_n})$ weakly star converges to a limit 
in $L^\infty(0,T;H)$ on one hand, and 
$$f_1'(\phi^{\alpha_n})\to f_1'(\phi)\quad\mbox{a.e.~in } Q$$
on the other. Then, the two limits must coincide and $f_1'(\phi)$ makes sense in $L^\infty(0,T;H)$,
although it is no longer expected that $f_1'(\phi)$ is uniformy bounded in $\LiQ$ (cf.~\eqref{ssboundal}).}

Therefore, passing to the limit as $n\to\infty$ in the state system \State\ written for $(\mu^{\alpha_n},\phi^{\alpha_n},w^{\alpha_n})$ and
$u^{\alpha_n}$, reveals that $(\mu,\phi,w)$ \pier{solves the system \eqref{ss01}--\eqref{ss05} at least in a variational sense. However, as $\mu \in \L\infty V$ and $ \dt\phi \in \L\infty H$, from \eqref{ss02}, \eqref{ss04} and standard elliptic regularity results if follows that $\mu \in \L\infty W$, whence 
$(\mu,\phi,w)$ is actually a strong solution of \eqref{ss01}--\eqref{ss05}.}
 According to \cite[Thm.~2.2]{CoSpTr}, this solution is uniquely determined, and thus 
we have $(\mu,\phi,w)=(\mu^0,\phi^0,w^0)$. Moreover, the unicity of the limit point entails that the  convergence properties are in 
fact valid for the entire sequence and not just for a subsequence.
This concludes the proof of the assertion.
\Edim

%\Brem
At this point, \pier{we note that} it makes sense to introduce the control-to-state operator for the viscous Cahn--Hilliard system 
\eqref{ss1}--\eqref{ss5}, 
\begin{equation}
\CS^0=(\CS^0_1,\CS^0_2,\CS^0_3):u^0\in\LiQ\mapsto (\CS^0_1(u^0),\CS^0_2(u^0),\CS^0_3(u^0)):=(\mu^0,\phi^0,w^0).
\end{equation}
%\Erem

\subsection{Convergence of optimal controls}
As we have seen \pier{in the previous subsection, the state variables associated with $\alpha\in (0,1]$ converge, as $\alpha\searrow0$ along a subsequence,} in a
well-defined sense to their counterparts for the case $\alpha=0$. We are now going to show that the same holds 
true for the optimal controls. 
\Bthm
\label{conv-opt-con}
Suppose that {\bf (A1)}--{\bf (A8)} are fulfilled, and let $\{\alpha_n\}_{n\in\enne}\subset (0,1]$ be given with $\alpha_n\searrow0$.
For every $n\in\enne$, let $u^{\alpha_n}\in\Uad$ be an optimal control for the problem {\bf (CP$_{\alpha_n}$)} with associated state
$(\mu^{\alpha_n},\phi^{\alpha_n},w^{\alpha_n})=\CS^{\alpha_n}(u^{\alpha_n})$. If, in addition, $u^{\alpha_n}\to u^0$ 
weakly star in $\LiQ$ as $n\to\infty$, then $u^0$ is an optimal control of the problem {\bf (CP$_0$)}.  
\Ethm
\Bdim
First we observe that, according to Theorem~\ref{connystate}, we have the convergence properties \pier{\eqref{connymu1}--\eqref{connyw}} where
$(\mu^0,\phi^0,w^0)=\CS^0(u^0)$, which means, in particular, that the pair $((\mu^0,\phi^0,w^0),u^0)$ is admissible for 
the limit problem {\bf (CP$_0$)}. Moreover, it obviously holds
\begin{align}
&\lim_{n\to\infty} \,\Big(\frac{b_1}2\iint_Q|\phi^{\alpha_n}-\phi_Q|^2\,+\,\frac{b_2}2\iO|\phi^{\alpha_n}(T)-\phi_\Omega|^2\Big) \non\\
&=\,\frac{b_1}2 \iint_Q |\phi^0-\phi_Q|^2\,+\,\frac{b_2}2\iO|\phi^0(T)-\phi_\Omega|^2\,.
\label{umba-umba}
\end{align}
Now let $u\in\Uad$ be arbitrary. \pier{We recall the weak sequential lower semicontinuity of the 
$L^2(Q)$-norm and of the functional $G$, which is convex and continuous on $L^2(Q)$. Then, 
from \eqref{umba-umba} and the optimality of the pair $((\mu^{\alpha_n},\phi^{\alpha_n},w^{\alpha_n}),
u^{\alpha_n})$ for the problem {\bf (CP$_{\alpha_n}$)}, it follows} that
\begin{align}
& \CJ(\CS^0_2(u^0),u^0)\,=\,\frac{b_1}2 \iint_Q |\phi^0-\phi_Q|^2\,+\,\frac{b_2}2\iO|\phi^0(T)-\phi_\Omega|^2\,+\,\frac{b_3}2\iint_Q|u^0|^2
\,+\,\kappa G(u^0)\non\\
&\le\,\liminf_{n\to\infty}\,\Bigl(\frac{b_1}2\iint_Q|\phi^{\alpha_n}-\phi_Q|^2\,+\,\frac{b_2}2\iO|\phi^{\alpha_n}(T)-\phi_\Omega|^2
\,+\,\frac{b_3}2\iint_Q|u^{\alpha_n}|^2 \,+\,\kappa G(u^{\alpha_n})\Bigr)\non\\[1mm] 
&=\,\liminf_{n\to\infty}\,\CJ(\CS^{\alpha_n}_2(u^{\alpha_n}),u^{\alpha_n})\,\le\,\liminf_{n\to\infty}\,\CJ(\CS^{\alpha_n}_2(u),u)
\,=\,\CJ(\CS^0_2(u),u)\,.
\end{align} 
Since $u\in\Uad$ was arbitrary, the assertion is proved.
\Edim
\Brem
1. Notice that any sequence $\{u^{\alpha_n}\}\subset\Uad$ with $\alpha_n\searrow 0$ contains at least one subsequence that converges in the
weak star topology of $\LiQ$ to some limit point $u\in\Uad$. From this it also follows that the problem {\bf (CP$_0$)} admits at least one
solution.  \\ 
2. Observe that we cannot expect a similar result for controls that are only locally optimal in the sense of $\LiQ$. Indeed, if $u^{\alpha_n}$ is 
locally optimal for {\bf (CP$_{\alpha_n}$)} in the sense of $\LiQ$, then we can only guarantee the existence of some $\varepsilon_n>0$ such 
that $\,\widehat \CJ(\CS^{\alpha_n}_2(u^{\alpha_n}),u^{\alpha_n})\,\le\, \widehat\CJ(\CS^{\alpha_n}_2(u),u)\,$ for all $u\in\Uad$ with
$\|u-u^{\alpha_n}\|_{\LiQ}\,\le\,\varepsilon_n$. Since it is well possible that $\,\varepsilon_n\to 0$ as $n\to\infty$, we
cannot conclude that limits of subsequences of $\{u^{\alpha_n}\}$ are locally optimal for {\bf (CP$_0$)}.    
\Erem

\subsection{Convergence of the adjoint variables}
Finally, we present a convergence result for the adjoint state variables \pier{under the following additional restriction.}

\vspace*{2mm}\noindent 
\pier{{\bf (A9)} \quad In the framework of {\bf (A3)} we further suppose that  $\, r_-= -\infty$ and $r_+= +\infty\,$ so that both $f_1$ and $f_2$  belong to $C^3(\erre)$.}

\Bthm
\label{conv-adj}
Suppose that {\bf (A1)}--\pier{{}{\bf (A9)}} are fulfilled, and let $\{\alpha_n\}_{n\in \enne}$ and $\{u^{\alpha_n}\}_{n\in\enne}\subset\Uad$ be sequences
with $\alpha_n\searrow0$ and $u^{\alpha_n}\to u^0$ weakly star in $\LiQ$. If $(\mu^{\alpha_n},\phi^{\alpha_n},w^{\alpha_n})
=\CS^{\alpha_n}(u^{\alpha_n})$ and $(p^{\alpha_n},q^{\alpha_n},r^{\alpha_n})$ denote the \pier{respective state and adjoint state variables associated}
with $u^{\alpha_n}$ for $n\in\enne$, then it holds
\begin{align}
\label{connyp1}
&p^{\alpha_n} \to p^0 \quad\mbox{weakly star in } \L\infty W,\\
\label{connyp2}
&\alpha_n\,p^{\alpha_n}\to 0 \quad\mbox{weakly star in } W^{2,\infty}(0,T;H) \mbox{ and strongly in } W^{1,\infty}(0,T;V),\\
\label{connyq}
&q^{\alpha_n}\to q^0 \quad\mbox{weakly star in } \L\infty V \cap\L2W,\\
\label{connypq} 
&p^{\alpha_n}+\tau q^{\alpha_n} \to p^0+\tau q^0 \quad \pier{\mbox{weakly star in } \,H^{1}(0,T;H)\cap \L\infty V\cap L^2(0,T;W)}
\non
\\
&\qquad\qquad\qquad\qquad\qquad \ \pier{\mbox{ and strongly in }\,\C0H \cap \L2V ,}\\
\label{connyr}
&r^{\alpha_n}\to r^0\quad\mbox{weakly star in } W^{1,\infty}(0,T;V)\cap \H1W,
\end{align}
where the triple $(p^0,q^0,r^0)$ is the unique solution to the adjoint system associated with the system \eqref{ss01}--\eqref{ss05} and the
control $u^0$, which is with $(\mu^0,\phi^0,w^0):=\CS^0(u^0)$ given by
\begin{align}
\label{as01}
&\pier{{}-\Delta p - q=0{}}
&&\mbox{a.e. in }Q\,,\\
\label{as02}
&\pier{{}-\dt(p+\tau q)-\Delta q+f''(\phi^0)q=b_1(\phi^0-\phi_Q){}}&&\mbox{a.e. in }Q\,,\\
\label{as03}
&-\gamma\dt r+r-q=0 &&\mbox{a.e. in }Q\,,\\
\label{as04}
&\dn p=\dn q =0 &&\mbox{a.e. on }\Sigma\,,\\
\label{as05}
& \pier{(p +\tau q) (T)={b_2}(\phi^0(T)-\phi_\Omega),  \,\,\, r(T)=0}&&\mbox{a.e. in }\Omega\,.
\end{align}
\Ethm
\Bdim
From Theorem~\ref{connystate} we infer that $(\mu^{\alpha_n},\phi^{\alpha_n},w^{\alpha_n})$ converges to $(\mu^0,\phi^0,w^0)$ in the sense
of \eqref{connymu1}--\eqref{connyw}. \pier{Moreover, recalling Remark~\ref{PIER1}, it turns out that under the validity of {\bf (A9)} the values $r_*$ and $ r^*$ 
in \eqref{separa} and the constant $K_2$ in \eqref{ssboundal} are independent of $\alpha$. Then, we can also infer that}
$$f''(\phi^{\alpha_n})\to f''(\phi^0) \quad\mbox{strongly in } C^0(\overline Q).$$
Moreover, it follows from \pier{the uniform estimate} \eqref{stimaAdj1} \pier{(see Remark~\ref{PIER2} as well)} that $(p^{\alpha_n},q^{\alpha_n},r^{\alpha_n})$ converges 
to a triple $(p,q,r)$ in the sense of 
\eqref{connyp1}--\eqref{connyr} at least on a subsequence. It is then readily seen that $(p,q,r)$ is a strong solution
to the system \eqref{as01}--\eqref{as05}.
In view of the linearity of this system, the solution is easily shown to be uniquely determined \pier{(anyway the reader may consult and refer to \cite[Theorem~4.3]{CoSpTr}),} and thus 
$(p,q,r)=(p^0,q^0,r^0)$. 

\pier{We point out the correctness of the first terminal  condition in \eqref{as05}
since the variable under time derivative in \eqref{as02} is just $(p+\tau q)$, so that the terminal condition follows from the convergence \eqref{connypq} and the conditions in \eqref{as5}.}

Hence, \pier{owing to the
unicity of the limit point $(p^0,q^0,r^0)$, the convergence properties are in fact valid for the entire sequence $\alpha_n\searrow0$} and not just for a subsequence.
This concludes the proof of the assertion. 
\Edim

%%%%%%%%%%%%%%%%%%%%%%%%%%%%%%%%%%%%%%%%%%%%%%%%%%%%%%%%%%%%%%%%%%%%%%%%%%%

\smallskip

\section*{Acknowledgments}
\pier{PC acknowledges the support of the Next Generation EU Project No.P2022Z7ZAJ 
(A unitary mathematical framework for modelling muscular 
dystrophies), the RISM (Research Institute for Mathematical Sciences,
an International Joint Usage/Research Center located in Kyoto  University)
and the GNAMPA (Gruppo Nazionale per l'Analisi
Matematica, la Probabilit\`{a} e le loro Applicazioni) of INdAM (Istituto
Nazionale di Alta Matematica).}
\medskip

\bigskip

%%%%%%%%%%%%%%%%%%%%%%%%%%%%%%%%%%%%%%%%%%%%%%%%%%%%%%%%%%%%%%%%%%%%%%%%%%%

\End{document}
%%%%%%%%%%%%%%%%%%%%%%%%%%%%%%%%%%%%%%%%%%%%%%%%%%%%%%%%%%%%%%%%%%%%%%%%%%%